\documentclass[draft]{amsart}

\usepackage{ifthen}

\usepackage[all]{xy}

\newtheorem{anyprop}{Anyprop}[section]

\newtheorem{theorem}[anyprop]{Theorem}
\newtheorem{lemma}[anyprop]{Lemma}

\newtheorem{corollary}[anyprop]{Corollary}

\theoremstyle{definition}

\newtheorem{example}[anyprop]{Example}

\newtheorem{remark}[anyprop]{Remark}

\newcommand{\ZZ}{\mathbb{Z}}

\newcommand{\FF}{\mathbb{F}}

\newcommand{\PP}{\mathbb{P}}

\newcommand  {\shE}     {\mathcal{E}}
\newcommand  {\shF}     {\mathcal{F}}

\newcommand  {\shL}     {\mathcal{L}}

\newcommand  {\shS}     {\mathcal{S}}

\newcommand  {\shX}     {\mathcal{X}}

\newcommand  {\fom}     {\mathfrak{m}}

\newcommand  {\fop}     {\mathfrak{p}}

\newcommand  {\Char}    {\operatorname{char}}

\newcommand  {\Ext}     {\operatorname{Ext}}

\renewcommand{\O}       {\mathcal{O}}

\newcommand  {\Pic}     {\operatorname{Pic}}

\newcommand  {\Proj}    {\operatorname{Proj}}

\newcommand  {\ra}      {\rightarrow}

\newcommand  {\rk}    {\operatorname{rk}}

\newcommand  {\Spec}    {\operatorname{Spec}}

\newcommand  {\Syz}     {\operatorname{Syz}}

\newcommand{\lto}{\longrightarrow}

\theoremstyle{remark}

\numberwithin{equation}{section}

\usepackage{amscd}
\usepackage{amssymb}

\begin{document}

\title[A remark on Frobenius descent for vector bundles]
{A remark on Frobenius descent for vector bundles}

\author[Holger Brenner and Almar Kaid]{Holger Brenner and Almar Kaid}

\address{Universit\"at Osnabr\"uck, Institut f\"ur Mathematik, 49069
Osnabr\"uck}

\address{Department of Pure Mathematics, University of Sheffield,
  Hicks Building, Hounsfield Road, Sheffield S3 7RH, United Kingdom}
\email{hbrenner@uni-osnabrueck.de and a.kaid@sheffield.ac.uk}




\begin{abstract}
We give a class of examples of vector bundles on a relative smooth
projective curve over $\Spec \ZZ$ such that for infinitely many
prime reductions the bundle has a Frobenius descent, but the
restriction to the generic fiber in characteristic zero is not
semistable. In the third section of the paper we prove for a large
class of varieties (including abelian varieties) that any vector
bundle with this Frobenius descent property is generically
semistable.
\end{abstract}

\maketitle



\section{Introduction}
Let $X$ be a smooth projective variety defined over an
algebraically closed field of characteristic $p>0$ with a fixed
very ample line bundle $\O_X(1)$. Further we denote by $F$ the
absolute Frobenius morphism $F:X \ra X$ which is the identity on
the topological space underlying $X$ and the $p$th power map on
the structure sheaf $\O_X$. A vector bundle $\shE$ on $X$ descends
under $F$ if there exists another vector bundle $\shF$ such that
$\shE \cong F^*(\shF)$. This note is inspired by the recent
preprint \cite{joshiremarkonvectorbundles} of K. Joshi. In the
relative situation, where a morphism $\shX \ra \Spec R$ with
generic fiber $X:=\shX_0$ is given and $R$ is a $\ZZ$-domain of
finite type, Joshi asked the following interesting question:
``assume $X$ is a smooth projective variety and suppose $V$ is a
vector bundle which descends under Frobenius modulo an infinite
set of primes then is it true that $V$ is semistable (with respect
to any ample line bundle on $X$)?'' He gives a positive answer to
this question for rank two vector bundles under the additional
assumption that $\Pic(X)=\ZZ$.

In section \ref{sectioncounterexample} of this paper we provide a
class of examples which give a negative answer to this question in
general. We show that on the relative Fermat curve
$C=V_+(X^d+Y^d+Z^d) \rightarrow \Spec \ZZ$, with $d \geq 5$ odd,
there exists a vector bundle $\shE$ of rank two such that for
infinitely many prime numbers $p$ the reduction $\shE_p=\shE|_{C_p}$
modulo $p$ has a Frobenius descent, but $\shE_0=\shE|_{C_0}$ is not
semistable on the fiber over the generic point.  In section
\ref{positiveresult} we give an affirmative answer to this
question under the assumption that for every closed point $\fom
\in \Spec R$ every semistable vector bundle on the fiber
$\shX_\fom$ is strongly semistable. We recall that a
semistable vector bundle $\shE$ is strongly semistable if
${F^e}^*(\shE)$ is semistable for all integers $e \geq 0$. This
provides further examples of varieties with $\Pic(X) \neq \ZZ$
(for example abelian varieties) for which the question of Joshi
still has a positive answer.

We would like to thank A. Werner for pointing out this problem to
us. Furthermore, we thank the referee for many useful comments
which helped to simplify the proof of Lemma \ref{hnfiltration} and
to clarify Example \ref{counterexample} via Lemma
\ref{notsplitlemma}.

\section{A counterexample for vector bundles on
curves}\label{sectioncounterexample}

In this section we give an example of a rank two vector bundle on
a generically smooth projective relative curve over $\Spec \ZZ$
such that infinitely many prime reductions have a Frobenius
descent but the bundle is not semistable on the generic fiber in
characteristic zero.

Our example will use the syzygy bundle $\Syz(X^2,Y^2,Z^2)(m)$ on
Fermat curves $C= V_+(X^d+Y^d+Z^d) \subset \PP^2$ defined over a
field $K$. This vector bundle is defined by the short exact
sequence
$$0 \lto \Syz(X^2,Y^2,Z^2)(m) \lto \O_C(m-2)^3 \lto \O_C(m) \lto
0,$$ where the penultimate mapping is given by $(s_1,s_2,s_3)
\mapsto s_1 X^2 +s_2 Y^2 + s_3 Z^2$. The bundle
$\Syz(X^2,Y^2,Z^2)(m)$ is semistable for $d \geq 5$ by
\cite[Proposition 6.2]{brennercomputationtight}. In positive
characteristic $p>0$, since the presenting sequence only involves
locally free sheaves, it is easy to see that the Frobenius
pull-back $F^*(\Syz(X^2,Y^2,Z^2)(m)) \cong
\Syz(X^{2p},Y^{2p},Z^{2p})(m p)$.

\begin{lemma}\label{hnfiltration}
Let $d=2\ell +1$ with $\ell \geq 2$ and let $C:=\Proj
K[X,Y,Z]/(X^d+Y^d+Z^d)$ be the Fermat curve of degree $d$ defined
over a field $K$ of characteristic $p \equiv \ell \mod d$. Then
the Frobenius pull-back of $\Syz(X^2,Y^2,Z^2)(3)$ sits inside the
short exact sequence
$$0 \lto \O_C(\ell-1) \lto \Syz(X^{2p},Y^{2p},Z^{2p})(3p) \lto
\O_C(-\ell+1) \lto 0.$$ In particular, the Frobenius pull-back is
not semistable and this sequence constitutes its Harder-Narasimhan
filtration.
\end{lemma}

\begin{proof}
We write $2p=dk+2\ell$ with $k$ even. The pull-back
$\Syz(X^{2p},Y^{2p},Z^{2p})$ of $\Syz(X^2,Y^2,Z^2)$ has a
non-trivial global section in total degree $d(k+1+k/2)$ by
\cite[Proof of Proposition 1.2]{brennermiyaoka}. From the
presenting sequence of the pull-back one reads off the degree as
follows:
\begin{eqnarray*}
\deg(\Syz(X^{2p},Y^{2p},Z^{2p})(d(k+1+k/2))&=& d(2d(k+1+k/2)-6p)\\
&=&d(2d(k+1+k/2)-3(dk+2\ell))\\
&=&d(2d-6\ell)\\
&=&d(-2\ell+2) <0.
\end{eqnarray*}
Since a semistable vector bundle of negative degree can not have
non-trivial global sections, the Frobenius pull-back
$\Syz(X^{2p},Y^{2p},Z^{2p})$ is not semistable. We obtain a
non-trivial mapping $\O_C(\ell-1) \ra
\Syz(X^{2p},Y^{2p},Z^{2p})(3p)$. We want to show that this mapping
constitutes the Harder-Narasimhan filtration of the pull-back,
meaning that this mapping has no zeros. Hence, assume that we have
a factorization
$$\O_C(\ell-1) \lto \shL \lto \Syz(X^{2p},Y^{2p},Z^{2p})(3p),$$
where $\shL$ is a subbundle of the syzygy bundle and has degree
$\deg(\shL):=\alpha \geq (\ell-1)d$. We have the short exact
sequence
$$0 \lto \shL \lto \Syz(X^{2p},Y^{2p},Z^{2p})(3p) \lto \shL^\prime \lto
0,$$ where $\shL^\prime$ is a line bundle of degree $-\alpha$. By
\cite[Corollary $2^p$]{shepherdbarronsemistability} (or
\cite[Theorem 3.1]{sunsemistable}) the inequality
$$\mu_{\max}(\shS) - \mu_{\min}(\shS) = \alpha - (-\alpha) = 2
\alpha \leq 2g-2$$  holds, where
$\shS:=\Syz(X^{2p},Y^{2p},Z^{2p})(3p)$ and $g$ denotes the genus
of $C$. The genus formula for plane curves yields
$$2g-2 =(d-1)(d-2)-2=d(d-3)=2d(\ell-1).$$
Therefore, we obtain $\alpha = d(\ell-1)$. Hence, $\O_C(\ell -1)
\cong \shL$ and the Harder-Narasimhan filtration is indeed $0
\subset \O_C(\ell-1) \subset \Syz(X^{2p},Y^{2p},Z^{2p})(3p).$
\end{proof}

\begin{remark}
Using Hilbert-Kunz theory and its geometric interpretation
developed in \cite{brennerhilbertkunz} and
\cite{trivedihilbertkunz} one can give an alternative (but more
complicated) proof that the line bundle $\O_C(\ell-1)$ is the
maximal destabilizing subbundle of the syzygy bundle
$\Syz(X^{2p},Y^{2p},Z^{2p})(3p)$. We recall that for a rank two
vector bundle the Harder-Narasimhan filtration is already strong
in the sense of \cite[Paragraph 2.6]{langersemistable}. By the
formula given in \cite[Theorem 3.6]{brennerhilbertkunz} we can
compute from the short exact sequence
$$0 \lto \shL \lto \Syz(X^{2p},Y^{2p},Z^{2p})(3p) \lto \shL^\prime \lto
0$$ the Hilbert-Kunz multiplicity $e_{HK}(I)$ (see
\cite{monskyhilbertkunz}) of the ideal $I=(X^2,Y^2,Z^2)$ in the
homogeneous coordinate ring $R:=K[X,Y,Z]/(X^d+Y^d+Z^d)$ of the
curve $C$ and obtain $e_{HK}(I)= 3d+\frac{\alpha^2}{dp^2}$. But,
by \cite[Theorem 2.3]{monskytrinomial} the Hilbert-Kunz
multiplicity of $I$ equals $e_{HK}(I)= 3d + \frac{d}{4}
\frac{(d-3)^2}{p^2}$ which implies $\alpha=d(\ell-1)$.
\end{remark}

\begin{remark}
We briefly comment on the situation for $\ell =0,1$. For $\ell =0$
(and $p \neq 2$) we have $\Syz(X^2,Y^2,Z^2)(3) \cong \O_{\PP^1}^2$
and this is also true for its Frobenius pull-back. For $\ell=1$,
we get the Fermat cubic which is an elliptic curve. In this case
we have an exact sequence $$0 \lto \O_C \lto \Syz(X^2,Y^2,Z^2)(3)
\lto \O_C \lto 0 \, ,$$ where the (only) global non-trivial
section is given by the curve equation. So the syzygy bundle is
$F_2$ in Atiyah's classification \cite{atiyahelliptic} and is
semistable, but not stable. Its Frobenius pull-back is either
$F_2$ (for $p \equiv 1 \mod 3$, i.e. Hasse invariant one) or
$\O_C^2$ (for $p \equiv 2 \mod 3$, i.e. Hasse invariant zero).
\end{remark}

In the relative situation
$$C:=\Proj(\ZZ_{d}[X,Y,Z]/(X^d+Y^d+Z^d)) \lto \Spec \ZZ_{d}$$
every fiber $C_p:=C \times_{\Spec \ZZ_{d}} \Spec \FF_p$ is a smooth
projective curve, namely the Fermat curve, defined over the prime
field $\FF_p$ (and $\bar{C}_p:=C \times_{\Spec \ZZ_{d}} \bar{\FF}_p$
is a smooth projective curve over the algebraic closure of
$\FF_p$) for every prime number $p$ such that $p \nmid d$. We
remind that by the Theorem of Dirichlet \cite[Chapitre VI, \S 4,
Th\'{e}or\`{e}me and Corollaire]{serrearithmetic} there exist
infinitely many prime numbers $p \equiv \ell \mod d$.

\begin{lemma}
\label{notsplitlemma} Let $d=2\ell+1$, $\ell \geq 2$, and consider
the smooth projective relative curve
$C:=\Proj(\ZZ_{d}[X,Y,Z]/(X^d+Y^d+Z^d)) \lto \Spec \ZZ_{d}$. Then the
sequence $($from Lemma \ref{hnfiltration}$)$ $$0 \lto
\O_{C_p}(\ell-1) \lto \Syz(X^{2p},Y^{2p},Z^{2p})(3p) \lto
\O_{C_p}(-\ell+1) \lto 0$$ does not split for almost all primes $p
\equiv \ell \mod d$.
\end{lemma}

\begin{proof}
Since $\Syz(X^{2p},Y^{2p},Z^{2p})(3p) \cong
F^*(\Syz(X^{2},Y^{2},Z^{2})(3))$ holds on every fiber $C_p$, the
bundle $\Syz(X^{2p},Y^{2p},Z^{2p})(3p)$ carries an integrable
connection $\nabla_p$ with $p$-curvature zero by the
Cartier-correspondence \cite[Theorem 5.1]{katznilpotent}. Assume
that the sequence does split for some $p \equiv \ell \mod d$. Then
$\O_{C_p}(\ell-1)$ is a direct summand of
$\Syz(X^{2p},Y^{2p},Z^{2p})(3p)$. The summand $\O_{C_p}(\ell-1)$
carries also a connection with the same properties. Hence, again
by the Cartier-correspondence it has a Frobenius descent and so
its degree $d(\ell-1)$ is divisible by $p$. But this can only hold
for finitely many $p$.
\end{proof}

\begin{example}
\label{counterexample}
As above we consider the smooth relative curve
$$C:= \Proj(\ZZ_{d}[X,Y,Z]/(X^d+Y^d+Z^d)) \lto \Spec \ZZ_{d},$$
with $d= 2 \ell +1$, $ \ell \geq 2$. The
$\rm\check{C}$ech-cohomology class $c=Z^{d-1}/XY \in
H^1(C,\O_C(d-3)) \cong \Ext^1(\O_C(-\ell+1),\O_C(\ell-1))$ defines
an extension
$$0 \lto \O_C(\ell-1) \lto \shE \lto \O_C(-\ell+1) \lto 0$$
with the corresponding restrictions to each fiber $C_\fop$, where
$\fop=(0)$ or $\fop=(p)$, $p \nmid d$. Note that this extension is
non-trivial on every fiber. This vector bundle $\shE$ is our
example. As $\ell \geq 2$ the bundle $\shE_0 = \shE|_{C_0}$ is not
semistable on $C_0$. By Lemma \ref{hnfiltration} we have for $p
\equiv \ell \mod d$ an extension
$$0 \lto \O_{C_p}(\ell-1) \lto \Syz(X^{2p},Y^{2p},Z^{2p})(3p) \lto
\O_{C_p}(-\ell+1) \lto 0$$ corresponding to $c^\prime \in
H^1(C_p,\O_{C_p}(2\ell-2))=H^1(C_p,\O_{C_p}(d-3))$, and by
Lemma \ref{notsplitlemma} we have $c^\prime \neq 0$ for almost all $p \equiv \ell \mod d$. We claim that $\shE_p = \shE|_{C_p} \cong F^*
(\Syz(X^2,Y^2,Z^2)(3))$ holds for these prime numbers. Since $\omega_{C_p}=\O_{C_p}(d-3)= \O_{C_p}(2 \ell -2)$
and $h^1(C_p,\omega_{C_p})=1$ we have $c = \lambda c^\prime$ for
some
$\lambda \in \FF_p^\times$. Moreover, multiplication with
$\lambda$ induces an automorphism $\omega_{C_p} \stackrel{\cdot
\lambda}{\to} \omega_{C_p}$ of line bundles as well as an
automorphism $H^1(C_p,\omega_{C_p}) \stackrel{\cdot \lambda} \ra
H^1(C_p,\omega_{C_p})$ of vector spaces.  We obtain a commutative
diagram
$$ \xymatrix {
0 \! \ar[r] & \O_{C_p}(2\ell-2) \! \ar[r] \ar[d]^{\cdot \lambda} &
\Syz(X^{2p},Y^{2p},Z^{2p})(3p +\ell -1)
\!\ar[r] \ar[d] &  \O_{C_p}\! \ar[r] \ar[d]^{=} & 0\\
0 \!\ar[r] & \O_{C_p}(2\ell-2) \!\ar[r] &   \shE_p(\ell -1) \!
\ar[r] & \O_{C_p}\!\ar[r] & 0 , }$$
where the map in the middle is an isomorphism of vector bundles.
Hence, $\shE_p \cong \Syz(X^{2p},Y^{2p},Z^{2p})(3p) \cong
F^*(\Syz(X^2,Y^2,Z^2)(3))$ and therefore $\shE_p$ admits a
Frobenius descent on every fiber $C_p$.
\end{example}

\begin{remark}
Example \ref{counterexample} extends to all Fermat curves $C^d =
V_+(X^d+Y^d+Z^d)$ where the degree $d$ has an odd divisor
$d^\prime \geq 5$. To see this we write $d=d^\prime n$ and look at
the cover $f: C^d \ra C^{d^\prime}$ induced by the ring map which
sends each variable to its $n$th power. Then the pull-back under
$f$ of the vector bundles considered in Example
\ref{counterexample} provide also an example on $C^d$ with the same properties.
\end{remark}

\section{A positive result}\label{positiveresult}

Let $\shX \ra \Spec R$ be a smooth projective morphism of relative
dimension $d\geq 1$, where $R$ is a domain of finite type over
$\ZZ$. Typical examples for the base are $\Spec \ZZ$ or arithmetic
schemes $\Spec D$, where $D$ is the ring of integers in a number
field. Let $\shE$ be a vector bundle over $\shX$. In \cite[Theorem
4.2]{joshiremarkonvectorbundles} K. Joshi proved under the
assumptions $\Pic(X)=\ZZ$ ($X=\shX_0$) and $\rk(\shE)=2$ that
$\shE_0=\shE|_{X}$ is semistable if for infinitely many closed
points $\fom \in \Spec R$ of arbitrarily large residue
characteristic the reduction $\shE_\fom$ admits a Frobenius
descent on the fiber $X_\fom= \shX_\fom$. The aim of this section
is to prove (using essentially the same methods) this result for
vector bundles of arbitrary rank under the assumption that for
every closed point $\fom$ every semistable vector bundle $\shF$ on
$X_\fom$ is strongly semistable, i.e. ${F^e}^*(\shF)$ is
semistable for all $e \geq 0$ (it is enough to assume this for
infinitely many closed points $\fom$ of arbitrary large residue
characteristic). It is interesting to note that Joshi used in
\cite[Theorem 2.1]{joshiremarkonvectorbundles} the condition
$\Pic(Y)=\ZZ$ on a smooth projective variety $Y$ in positive
characteristic and a further hypothesis on $Y$ to prove that every
semistable rank two vector bundle on $Y$ is strongly semistable.

\begin{theorem}
\label{frobeniusdescentcriterion} Let $R$ be a $\ZZ$-domain of
finite type and let $f: \shX \rightarrow \Spec R$ be a smooth
projective morphism of relative dimension $d \geq 1$ together with a
fixed $f$-very ample line bundle $\O_\shX(1)$ and let $\shE$ be a
vector bundle on $\shX$. Further assume that every semistable vector
bundle is strongly semistable $($with respect to
$\O_{X_{\fom}}(1)$$)$ for every fiber $X_{\fom}$, $\fom$ a closed
point in $\Spec R$. Then the following holds: If $\shE_\fom =
\shE|_{X_\fom}$ has a Frobenius descent for infinitely many closed
points $\fom \in \Spec R$ of arbitrarily large residue
characteristic, then $\shE_0$ is semistable on the generic fiber
$X=X_0=\shX_0$.
\end{theorem}

\begin{proof}
One can show by induction over $\dim R$ that there exists a bound
$b$ such that $\mu_{\max}(\shE_\fom) \leq b$ for all closed points
$\fom \in \Spec R$ (see \cite[Lemma
3.1]{brennerkaiddeepfrobeniusdescent} for an explicit proof). For a
closed point $\fom \in \Spec R$ with descent data $\shE_\fom \cong
F^*(\shF_\fom)$, $\shF_\fom$ locally free on the fiber $X_\fom$, we
have $$\mu_{\max}(\shE_\fom)=\Char(\kappa(\fom))
\mu_{\max}(\shF_\fom)$$ because semistable vector bundles are
strongly semistable on every fiber $X_\fom$ by assumption. Since we
have $\shE_{\fom} \cong F^*(\shF_{\fom})$ for infinitely many closed
points $\fom$ of arbitrarily large residue characteristics, this
forces the similar equalities $\deg(\shE_0) =\deg(\shE_\fom) =
\Char(\kappa(\fom))
\deg(\shF_\fom)$ (we take the degree always with respect to
$\O_{X_\fom}(1)$) which implies $\deg(\shE_\fom)=\deg(\shF_\fom)=0$.
Assume the restriction $\shE_0$ to the generic fiber $X$ is not
semistable. Then by the openness of semistability \cite[Section
5]{miyaokachern} every restriction
$\shE_\fom$ on $X_\fom$ is not semistable. Again by our assumption,
$\shF_\fom$ is not semistable either and so $\mu_{\max}(\shF_\fom)
\geq 1/r$, $r= \rk(\shE)$. This gives
$$b \geq \mu_{\max}(\shE_\fom)=\Char(\kappa(\fom))
\mu_{\max}(\shF_\fom) \geq \frac{\Char(\kappa(\fom))}{r}$$ which
contradicts the assumption that we have Frobenius descent at closed
points $\fom \in \Spec R$ of arbitrarily large residue
characteristic.
\end{proof}

\begin{corollary}
Let $R$ be a $\ZZ$-domain of finite type and let $f: \shX
\rightarrow \Spec R$ be a smooth projective morphism of relative
dimension $d \geq 1$ together with a fixed $f$-very ample line
bundle $\O_\shX(1)$ and let $\shE$ be a vector bundle on $\shX$.
Suppose that the fibers $X_\fom$, $\fom \in \Spec R$ closed,
fulfill at least one of the following $($not necessarily
independent$)$ properties:
\begin{enumerate}
\item $X_\fom$ is an abelian variety, \item $X_\fom$ is a
homogenous space of the form $G/P$ where $P$ is a reduced
parabolic subgroup, \item the cotangent bundle $\Omega_{X_\fom}$
fulfills $\mu_{\max}(\Omega_{X_{\fom}}) \leq 0$.
\end{enumerate}
Then the following holds: If $\shE_\fom$ has a Frobenius descent
for infinitely many closed points $\fom \in \Spec R$ of
arbitrarily large residue characteristics, then $\shE_0$ is
semistable on $X=X_0$.
\end{corollary}

\begin{proof}
That every semistable vector bundle is strongly semistable in the
case $(3)$ is due to \cite[Theorem
2.1]{mehtaramanathanhomogeneous}. Condition $(3)$ holds in
particular for the varieties occurring in $(1)$ and $(2)$. Other
proofs for this property in cases $(1)$ and $(2)$ are given in
\cite[Corollary $3^p$]{shepherdbarronsemistability} and in case
$(3)$ in \cite[Corollary 6.3]{langersemistable}. Hence, the
assertion follows from Theorem \ref{frobeniusdescentcriterion}.
\end{proof}

\begin{remark}
On the one hand it is well-known that every semistable vector
bundle on an elliptic curve is strongly semistable (cf.
\cite[Appendix]{tusemistable}). So elliptic curves provide an
important class of smooth projective varieties with $\Pic(X)\neq
\ZZ$ for which Theorem \ref{frobeniusdescentcriterion} holds. On
the other hand it is also known that for every smooth projective
curve of genus $g \geq 2$ there exists a semistable vector bundle
$\shF$ such that $F^*(\shF)$ is not semistable (see \cite[Theorem
1]{langepaulyfrobenius}). So we see that Theorem
\ref{frobeniusdescentcriterion} is applicable in relative
dimension one only for elliptic curves and the projective line
$\PP^1$.
\end{remark}

\bibliographystyle{amsplain}

\end{document}